\documentclass[12pt]{article}
\usepackage{amssymb}

\usepackage{amsmath,latexsym,amsfonts}
\usepackage{graphicx}
\usepackage{srcltx}
\usepackage{amscd}
\newtheorem{lemma}{Lemma}[section]
\newtheorem{coro}[lemma]{Corollary}
\newtheorem{prop}[lemma]{Proposition}
\newtheorem{thm}[lemma]{Theorem}
\newtheorem{defn}[lemma]{Definition}
\makeatletter\@addtoreset{equation}{section}
\renewcommand\theequation{\thesection.\@arabic\c@equation}

\begin{document}
\begin{center}
{\LARGE On transversally
harmonic maps of foliated Riemannian manifolds}

 \renewcommand{\thefootnote}{}
\footnote{${}^*$ Corresponding author}
\renewcommand{\thefootnote}{\arabic{footnote}}
\setcounter{footnote}{0}

\vspace{0.5 cm} {\large Min Joo Jung and Seoung Dal Jung${}^*$}
\end{center}
\vspace{0.5cm}
\noindent {\bf Abstract.} Let $(M,\mathcal F)$ and $(M',\mathcal F')$ be two foliated Riemannian manifolds with $M$ compact. Then we give a new proof of the first normal variational formula for the trasnversal energy.  Moreover, if we assume that the transversal Ricci curvature of $\mathcal F$ is nonnegative and the transversal sectional curvature of $\mathcal F'$
is nonpositive, then  any transversally harmonic map 
$\phi:(M, \mathcal
F) \rightarrow (M', \mathcal F')$ is  transversally totally geodesic. In addition, if the transversal Ricci curvature is positive at some point, then $\phi$ is trasnversally constant. 
\section{Introduction}
Transversally harmonic maps of foliated Riemannian manifolds were introduced by Konderak and Wolak [\ref{KW1}] in 2003. 
Let $(M,\mathcal F)$ and $(M',\mathcal F')$ be two foliated Riemannian manifolds and let $\phi:M\to M'$ be a smooth foliated map, i.e., $\phi$ is a smooth leaf-preserving map. Then $\phi$ is said to be {\it transversally harmonic} if the transversal tension field $\tau_b(\phi)$ vanishes. See Section 3 and [\ref{KW1}] for details. Equivalently, it is a critical point of the transversal energy functional on any compact domain of $M$, which is defined in Section 4 (cf.[\ref{KW2}]). Also, transversally harmonic maps are considered as harmonic maps between the leaf spaces [\ref{KW1},\ref{KW2}]. So, for the point foliation,  transversally harmonic maps are  harmonic maps. Therefore transversally harmonic maps are considered as generalizations of harmonic maps. In this paper, we study transversally harmonic maps and give some interesting facts relating to them. The paper is organized as follows.  In Section 2, we review the well-known facts on a foliated Riemannian manifold. In Section 3, we review the properties of the transversally harmonic map, which were studied in [\ref{KW2}] and give some results. In Section 4, we give a new proof of the first normal variational formula for the transversal energy $E_B(\phi)$ (Theorem 4.2). In the last section, we study the generalized Weitzenb\"ock formula and give some applications (Theorem 5.3 and Theorem 5.4).

\section{Preliminaries}

Let $(M,g,\mathcal F)$ be a $(p+q)$-dimensional  foliated Riemannian
manifold with a foliation $\mathcal F$ of codimension $q$ and a bundle-like metric $g$ with respect to $\mathcal F$ [\ref{Molino}, \ref{Tond}]. A {\it foliated Riemannian manifold}  means a Riemannian
manifold with a Riemannian foliation. Let $TM$ be the tangent bundle of $M$, $L$
the tangent bundle of $\mathcal F$, and  $Q=TM/L$ the
corresponding normal bundle of $\mathcal F$.
 Then we have an exact sequence of vector bundles
\begin{align}\label{eq1-1}
 0 \longrightarrow L \longrightarrow
TM_{\buildrel \longleftarrow \over \sigma }^{\buildrel \pi \over
\longrightarrow} Q \longrightarrow 0,
\end{align}
where $\pi:TM\to Q$ is a projection and $\sigma:Q\to L^\perp$ is a bundle map satisfying
$\pi\circ\sigma=id$. Let $g_Q$ be the holonomy invariant metric on
$Q$ induced by $g=g_L + g_{L^\perp}$; that is,
\begin{equation}
g_Q(s,t)=g(\sigma(s),\sigma(t)) \quad \forall\ s,t\in \Gamma Q.
\end{equation}
This means that $\theta(X)g_Q=0$ for $X\in\Gamma L$, where
$\theta(X)$ is the transverse Lie derivative. So we have an
identification $L^\perp$ with $Q$ via an isometric splitting
$(Q,g_Q)\cong (L^\perp, g_{L^\perp})$. We denote by $\nabla^Q$ the transverse Levi-Civita 
connection on the normal bundle $Q$ [\ref{Tond},\ref{Tond1}]. The transversal curvature tensor $R^Q$ of $\nabla$ is defined by $R^Q(X,Y)=[\nabla_X,\nabla_Y]-\nabla_{[X,Y]}$ for any $X,Y\in\Gamma TM$. It is trivial that $i(X)R^Q=0$ for any $X\in\Gamma L$. Let $K^Q,{\rm Ric}^Q $ and
$\sigma^Q$ be the transversal
sectional curvature, transversal Ricci operator and  transversal
scalar curvature with respect to $\nabla^Q\equiv\nabla$, respectively. The foliation $\mathcal F$ is said to be {\it
minimal} if $\kappa=0$, where $\kappa$ is the mean curvature form of $\mathcal F$ [\ref{Tond}].
Let $\Omega_B^r(\mathcal F)$ be the space of all {\it basic
$r$-forms}, i.e.,  $\phi\in\Omega_B^r(\mathcal F)$ if and only if
$i(X)\phi=0$ and $\theta(X)\phi=0$ for any $X\in\Gamma L$, where
$i(X)$ is the interior product. Then $\Omega^*(M)=\Omega_B^*(\mathcal F)\oplus \Omega_B^*(\mathcal F)^\perp$ [\ref{Lop}].  Let $\kappa_B$ be the basic part of $\kappa$. Then $\kappa_B$ is closed, i.e., $d\kappa_B=0$ [\ref{Lop}]. Now,  we define  the {\it basic
Laplacian} $\Delta_B$ acting on $\Omega_B^*(\mathcal F)$ by
\begin{equation}\label{eq1-5}
\Delta_B=d_B\delta_B+\delta_B d_B,
\end{equation}
where $\delta_B$ is the formal adjoint of
$d_B=d|_{\Omega_B^*(\mathcal F)}$ [\ref{Park}]. Let $\{E_a\}_{a=1,\cdots,q}$ be a local orthonormal basic frame on $Q$. We define $\nabla_{\rm tr}^*\nabla_{\rm tr}:\Omega_B^r(\mathcal F)\to
\Omega_B^r(\mathcal F)$ by
\begin{align}\label{eq1-7}
\nabla_{\rm tr}^*\nabla_{\rm tr} =-\sum_a \nabla^2_{E_a,E_a}
+\nabla_{\kappa_B^\sharp},
\end{align}
where $\nabla^2_{X,Y}=\nabla_X\nabla_Y -\nabla_{\nabla^M_XY}$ for
any $X,Y\in TM$. Then the operator $\nabla_{\rm tr}^*\nabla_{\rm tr}$
is positive definite and formally self adjoint on the space of
basic forms [\ref{Jung}].
 Let $V(\mathcal F)$ be the space of all transversal infinitesimal automorphisms $Y$ of $\mathcal F$, i.e., $[Y,Z]\in \Gamma L$ for
all $Z\in \Gamma L$ [\ref{Kamber2}]. Let
\begin{equation}
\bar V(\mathcal F)=\{\bar Y=\pi(Y)| Y\in V(\mathcal F)\}.
\end{equation}
Note that $\bar V(\mathcal F)\cong \Omega_B^1(\mathcal F)$[\ref{PP}].
For later use, we recall the transversal divergence theorem
[\ref{Yorozu}] on a foliated Riemannian
manifold.
\begin{thm} \label{thm1-1} $(${\rm Transversal divergence theorem}$)$
Let $(M,g_M,\mathcal F)$ be a closed, oriented Riemannian manifold
with a transversally oriented foliation $\mathcal F$ and a
bundle-like metric $g_M$ with respect to $\mathcal F$. Then
\begin{equation}
\int_M \operatorname{div_\nabla}\bar X = \int_M g_Q(\bar
X,\kappa_B^\sharp)
\end{equation}
for all $X\in V(\mathcal F)$, where $\operatorname{div_\nabla} X$
denotes the transversal divergence of $X$ with respect to the
connection $\nabla^Q$.
\end{thm}
Now we define the bundle map $A_Y:\Lambda^r
Q^*\to\Lambda^r Q^*$ for any $Y\in V(\mathcal F)$ [\ref{Kamber2}]
by
\begin{align}\label{eq1-11}
A_Y\phi =\theta(Y)\phi-\nabla_Y\phi.
\end{align}
It is well-known [\ref{Kamber2}] that for any $s\in\Gamma Q$
\begin{align}\label{2-5-1}
A_Y s = -\nabla_{Y_s}\bar Y,
\end{align}
where $Y_s$ is the vector field such that $\pi(Y_s)=s$. So $A_Y$ depends only on $\bar Y=\pi(Y)$.
Since
$\theta(X)\phi=\nabla_X\phi$ for any $X\in\Gamma L$, $A_Y$
preserves the basic forms and depends only on $\bar Y$. Now, we
recall the generalized Weitzenb\"ock formula on $\Omega_B^*(\mathcal F)$.
\begin{thm} $[\ref{Jung}]$ On a foliated Riemannian manifold $(M,\mathcal F)$, we have
\begin{align}\label{eq1-12}
  \Delta_B \phi = \nabla_{\rm tr}^*\nabla_{\rm tr}\phi +
  F(\phi)+A_{\kappa_B^\sharp}\phi,\quad\phi\in\Omega_B^r(\mathcal
  F),
\end{align}
 where $F(\phi)=\sum_{a,b}\theta^a \wedge i(E_b)R^\nabla(E_b,
 E_a)\phi$. If $\phi$ is a basic 1-form, then $F(\phi)^\sharp
 ={\rm Ric}^Q(\phi^\sharp)$.
\end{thm}
Now we recall a very important lemma for later use. From Proposition  4.1 in [\ref{Park}], it is well-known that $\Delta_B-\kappa_B^\sharp$ on all basic functions is the restriction of $\Delta-\kappa^\sharp$ on all functions. Hence, by maximum and minimum principles, we have the following lemma.
\begin{lemma}\label{lem2-3} $[\ref{JLK}]$ Let $(M,g,\mathcal F)$ be a closed, connected Riemannian manifold with a foliation $\mathcal F$ and a bundle-like metric $g$. If $(\Delta_B-\kappa_B^\sharp)f\geq 0$ $($or $\leq 0)$ for any basic function $f$, then $f$ is constant.
\end{lemma}

\section{Transversally harmonic maps}
Let $(M,  g,\mathcal F)$  and $(M', g',\mathcal F')$  be
two foliated Riemannian manifolds.
  Let $\nabla^{M}$ and $\nabla^{M'}$ be the Levi-Civita connections of
$M$ and $M'$, respectively. Let $\nabla$ and $\nabla'$ be the transverse
Levi-Civita connections on $Q$ and $Q'$, respectively.  Let $\phi:(M,g,\mathcal
F)\to (M', g',\mathcal F')$ be a smooth {\it foliated map},
i.e., $d\phi(L)\subset L'$.  Then  we define $d_T\phi:Q \to Q'$ by
\begin{align}
d_T\phi := \pi' \circ d \phi \circ \sigma. 
\end{align}
Then $d_T\phi$ is a section in $ Q^*\otimes
\phi^{-1}Q'$, where $\phi^{-1}Q'$ is the pull-back bundle on $M$. Let $\nabla^\phi$
and $\tilde \nabla$ be the connections on $\phi^{-1}Q'$ and
$Q^*\otimes \phi^{-1}Q'$, respectively. Then a foliated map $\phi:(M,\mathcal F)\to (M',\mathcal F')$ is called {\it transversally totally geodesic} if it satisfies
\begin{align}
\tilde\nabla_{\rm tr}d_T\phi=0,
\end{align}
where $(\tilde\nabla_{\rm tr}d_T\phi)(X,Y)=(\tilde\nabla_X d_T\phi)(Y)$ for any $X,Y\in \Gamma Q$. Note that if $\phi:M\to M'$ is transversally totally geodesic with $d\phi(Q)\subset Q'$, then, for any transversal geodesic $\gamma$ in $M$, $\phi\circ\gamma$ is also transversal geodesic. The {\it
transversal tension field} of $\phi$ is defined by
\begin{align}\label{eq3-2}
\tau_b(\phi)={\rm tr}_{Q}\tilde \nabla d_T
\phi=\sum_{a=1}^{q}(\tilde\nabla_{E_a} d_T\phi)(E_a),
\end{align}
where $\{E_a\}$ is a local orthonormal basic frame of $Q$.
Trivially, the transversal tension field $\tau_b(\phi)$ is a
section of $\phi^{-1}Q'$.
\begin{defn}{\rm
Let $\phi: (M, g,\mathcal F) \to (M', g',\mathcal F')$ be
a smooth foliated map. Then  $\phi$ is said to be }
transversally harmonic {\rm if the transversal tension field of $\phi$
vanishes}, i.e., $\tau_b(\phi)=0$. 
\end{defn} 
Now we recall the O'Neill tensors $\cal A$ and $\cal T$ [\ref{Oneill}] on a foliated manifold $(M,\mathcal F)$, which are
 defined by
\begin{align}
{\cal A}_XY &= \pi^\perp (\nabla_{\pi(X)}^M \pi(Y)) + \pi
(\nabla_{\pi(X)}^M \pi^\perp(Y))\label{eq2-7} \\
{\cal T}_XY &= \pi^\perp (\nabla_{\pi^\perp (X)}^M \pi(Y)) + \pi
(\nabla_{\pi^\perp(X)}^M \pi^\perp(Y))\label{eq2-8}
\end{align}
for any $X,Y\in TM$,
where $\pi^\perp : TM \to L$. It is well-known [\ref{Oneill}] that
\begin{align}\label{eq2-9}
{\cal A}_{\pi(X)}\pi(Y)=\pi^\perp[\pi(X),\pi(Y)]
\end{align}
for any vector fields $X,Y$ on $M$. 
  Then ${\cal T}\equiv 0$
is equivalent to the property that all leaves of $\mathcal F$ are
totally geodesic submanifolds of $(M, g)$ and ${\cal A}\equiv 0$ is
equivalent to the integrability of $Q$. 

    Let $\{E_i\}_{i=1,\cdots,p}$ be a local orthonormal basis of $L$ and 
   $\{E_a\}_{a=1,\cdots,q}$ be a local orthonormal basic frame on $Q$. Then we have the following.
\begin{thm}
Let $\phi:(M, g, \mathcal F) \rightarrow (M', g', \mathcal
F')$ be a smooth foliated map. Then
\begin{align*}
\tau(\phi) &= \tau(\phi|_{\mathcal
F})+\tau_b(\phi)-d_T\phi(\kappa^\sharp)+{\rm tr}_{g}\phi^* \mathcal T'+{\rm tr}_{Q}\phi^*\mathcal A'\\
&  + \sum_a\{\pi^\perp \nabla_{\pi d\phi(E_a)}^{M'} \pi^\perp
d\phi(E_a) +\pi^\perp \nabla_{\pi^\perp d\phi(E_a)}^{M'}
\pi^\perp d\phi(E_a) - \pi^\perp d\phi(\nabla_{E_a}
E_a)\}\\
&+ \sum_a \pi\nabla^{M'}_{\pi^\perp d\phi(E_a)}\pi d\phi(E_a),
\end{align*}
where $\tau(\phi)$ is the tension field of $\phi$ and 
\begin{align}
\tau(\phi|_{\mathcal F})=\pi^\perp\sum_i (\tilde\nabla_{E_i}d\phi)(E_i).
\end{align}
\end{thm}
{\bf Proof.} Let $\{E_i, E_a \}_{i=1,\cdots,p; a=1,\cdots,q}$ be a local
orthonormal frame of $TM$ such that $E_i\in \Gamma L$, $E_a\in\Gamma Q$. By the definition of the tension field, we have
\begin{align}\label{eq3-8}
\tau(\phi)&=\sum_{i=1}^{p}(\tilde\nabla_{E_i} d\phi)(E_i) + \sum_{a=1}^{q}(\tilde\nabla_{E_a} d\phi)
(E_a).
\end{align}
 Since $\phi$ is a foliated map, $\pi d\phi (E_i)=0$ and $\pi^\perp d\phi(E_i)=d\phi(E_i)$.
Therefore, we have
\begin{align*}
\sum_{i=1}^{p}(\tilde\nabla_{E_i} d\phi)(E_i)&=\tau(\phi|_{\mathcal F}) +\sum_i\{\pi\nabla_{d\phi(E_i)}^{M'} d\phi(E_i)
-\pi d\phi(\nabla_{E_i}^{M} E_i)\}
\end{align*}
and
\begin{align*}
&\sum_{a=1}^{q}(\tilde\nabla_{E_a} d\phi)
(E_a)\\&=\tau_b(\phi)+\sum_a\{\pi^\perp\nabla^{M'}_{\pi d\phi(E_a)}\pi d\phi(E_a) +\nabla^{M'}_{\pi d\phi(E_a)}\pi^\perp d\phi(E_a)\} \\ 
&  +\sum_a\{\nabla_{\pi^\perp d\phi(E_a)}^{M'} \pi d\phi(E_a)+
\nabla_{\pi^\perp d\phi(E_a)}^{M'} \pi^\perp d\phi(E_a)
-\pi^\perp d\phi(\nabla_{E_a}^{M}E_a)\}.
\end{align*}
From (\ref{eq2-9}), we have $\pi^\perp\nabla^{M'}_{\pi d\phi(E_a)}\pi d\phi(E_a)=\pi^\perp\nabla^{M}_{E_a}E_a=0$. Hence, from (\ref{eq2-7}) and (\ref{eq2-8}), we have
\begin{align*}
\tau(\phi)&=\tau(\phi|_{\mathcal F})+\tau_b(\phi)-\pi d\phi(\sum_i\pi(\nabla^{M}_{E_i}E_i))+
\sum_i {{\cal T}'}_{d\phi(E_i)}d\phi(E_i)\\
& +\sum_a \{{{\cal T}'}_{d\phi(E_a)}d\phi(E_a)+{{\cal A}'}_{d\phi(E_a)}d\phi(E_a)+ \pi\nabla^{M'}_{\pi^\perp d\phi(E_a)}\pi d\phi(E_a)\}\\
&+\sum_a\pi^\perp\{\nabla^{M'}_{\pi^\perp d\phi(E_a)}{\pi^\perp d\phi(E_a)}+\nabla^{M'}_{\pi d\phi(E_a)}{\pi^\perp d\phi(E_a)}-d\phi(\pi\nabla^{M}_{E_a}E_a)\}.
\end{align*}
Since $\sum_i \pi(\nabla^{M}_{E_i}E_i) =\kappa^\sharp$, the proof is completed. $\Box$
\begin{coro} $[\ref{X}]$ Let $\phi:(M, g, \mathcal F) \rightarrow (M', g', \mathcal
F')$ be a smooth foliated map and $d\phi(Q)\subset Q'$. Then
\begin{align*}
\tau(\phi) = \tau(\phi|_{\mathcal F}) + \tau_b (\phi) -d\phi(\kappa^\sharp) +
{\rm tr}_{L} \phi^*{\mathcal T'}.
\end{align*}
\end{coro}
{\bf Proof.} Since $d\phi(Q)\subset Q'$, $\pi^\perp d\phi(E_a)=0$ for all $a$. Moreover, from (\ref{eq2-8}) and (\ref{eq2-9}), ${\mathcal A'}_X X=0$ and ${\mathcal T'}_X Y=0$ for any $X,Y\in Q'$. Hence the proof is completed. $\Box$
\begin{coro} $[\ref{KW2}]$
Let $\phi:(M, g, \mathcal F) \rightarrow (M', g', \mathcal
F')$ be a  smooth foliated map and $d\phi(Q)\subset Q'$. Assume that $\mathcal F$ is
minimal and $\mathcal F'$ is totally geodesic. Then $\phi$ is harmonic if and only if $\phi$ is transversally harmonic and leaf-wise harmonic, i.e., $\tau(\phi|_{\mathcal F})=0$.
\end{coro}
{\bf Proof.} Since $\mathcal F$ is minimal and $\mathcal F'$ is totally geodesic, i.e, $\kappa=0$ and $\mathcal T'=0$, by Corollary 3.3, we have
\begin{align*}
\tau(\phi)=\tau(\phi|_{\mathcal F})+\tau_b(\phi).
\end{align*}
So the proof is completed. $\Box$
\begin{coro}
Let $\phi:(M, g, \mathcal F) \rightarrow (M', g', \mathcal
F')$ be a smooth foliated map and $d\phi(Q)\subset Q'$. Then $\phi$ is a transversally
harmonic map if and only if 
\begin{align*}
\pi(\tau(\phi))={\rm tr}_{L}\phi^*\mathcal T'
-d\phi(\kappa^\sharp).
\end{align*}
\end{coro}
Now, let $\mathcal F$ be a Riemannian flow defined by a unit vector field $V$ on a Riemannian manifold $(M^{n+1},g)$. Then
\begin{align}
\kappa^\sharp =\pi(\nabla^M_VV)=\nabla^M_VV.
\end{align}
In fact, $\nabla^M_VV$ is already orthogonal to the leaves since $g(\nabla^M_VV,V)=0$. Moreover, it is trivial that 
$\mathcal F$ is totally geodesic if and only if $\mathcal F$ is minimal, i.e., $\mathcal T=0$ if and only if $\kappa^\sharp=0$.    
    Let $\mathcal F$ and $\mathcal F'$ be two Riemannian flows defined by unit vector fields $V$ and $V'$ on Riemannian manifolds $(M,g)$ and $(M',g')$, respectively. Let $\phi:(M,\mathcal F)\to (M',\mathcal F')$ be a smooth foliated map. Then
\begin{align}
\tau(\phi|_{\mathcal F}) = V(\lambda)V' -\pi^\perp d\phi(\kappa^\sharp),\quad\lambda =(\phi^*\omega')(V),
\end{align}
where $\omega'$ is a dual 1-form of $V'$.
Hence $\phi$ is leaf-wise harmonic and $d\phi(Q)\subset Q'$ if and only if $\lambda$ is basic, i.e., $V(\lambda)=0$.   
Hence we have the following corollary.
\begin{coro} Let $\mathcal F$ and $\mathcal F'$ be two Riemannian flows defined by a unit vector field $V$ and $V'$ on a Riemannian manifold $M$ and $M'$, respectively. Assume that $\mathcal F$ and $\mathcal F'$ are minimal. Let $\phi:(M,g,\mathcal F)\to (M',g',\mathcal F')$ be a smooth foliated map and $d\phi(Q)\subset Q'$. Then $\phi$ is harmonic if and only if $\phi$ is transversally harmonic and $\phi^*(\omega')(V)$ is basic.
\end{coro}
Let $\phi:(M,\mathcal F)\to (M',\mathcal F')$ and $\psi:(M',\mathcal F')\to (M'',\mathcal F'')$ be smooth foliated maps. Then the composition $\psi\circ \phi:(M,\mathcal F)\to (M'',\mathcal F'')$ is a smooth foliated map. Moreover, we have
\begin{align}\label{eq3-11}
d_T (\psi\circ\phi)=d_T\psi \circ d_T\phi.
\end{align}
Hence we have the following proposition.
\begin{prop} Let $\phi:(M,\mathcal F)\to (M',\mathcal F')$ and $\psi:(M',\mathcal F')\to (M'',\mathcal F'')$ be smooth foliated maps. Then 
\begin{align}\label{eq3-12}
\tilde\nabla_{\rm tr}d_T (\psi\circ\phi)=d_T\psi (\tilde\nabla_{\rm tr}d_T\phi) +\phi^*\tilde\nabla_{\rm tr}d_T\psi,
\end{align}
where $(\phi^*\tilde\nabla_{\rm tr}d_T\psi)(X,Y)=(\tilde\nabla_{d_T\phi(X)}d_T\psi )(d_T\phi(Y))$ for any $X,Y\in \Gamma Q$.
\end{prop}
{\bf Proof.} From (\ref{eq3-11}),  we have that, for any  $X,Y \in\Gamma Q$,
\begin{align*}
(\tilde\nabla_{\rm tr} d_T(\psi\circ\phi))(X,Y)&=\nabla_X^{\psi\circ\phi}d_T(\psi\circ\phi)(Y)-d_T(\psi\circ\phi)(\nabla_X Y)\\
&=(\tilde\nabla_{d_T\phi(X)}d_T\psi)(d_T\phi(Y))+d_T\psi((\tilde\nabla_X d_T\phi)(Y))\\
&=(\phi^*\tilde\nabla_{\rm tr}d_T\psi)(X,Y) +d_T\psi(\tilde\nabla_{\rm tr}d_T\phi)(X,Y),
\end{align*}
which proves (\ref{eq3-12}). $\Box$
\begin{coro} Let $\phi:(M,\mathcal F)\to (M',\mathcal F')$ and $\psi:(M',\mathcal F')\to (M'',\mathcal F'')$ be smooth foliated maps. Then the transversal tension field of the composition is given by
\begin{align}
\tau_b(\psi\circ\phi)=d_T\psi(\tau_b(\phi))+{\rm tr}_Q \phi^*\tilde\nabla_{\rm tr}d_T\psi,
\end{align}
where ${\rm tr}_Q \phi^*\tilde\nabla_{\rm tr}d_T\psi =\sum_{a=1}^q(\tilde\nabla_{d_T\phi(E_a)}d_T\psi)(d_T\phi(E_a))$.
\end{coro}
\begin{coro} Let $\phi:(M,\mathcal F)\to (M',\mathcal F')$ be a transversally harmonic map and let $\psi:(M',\mathcal F')\to (M'',\mathcal F'')$ be a transversally totally geodesic map. Then $\psi\circ\phi:(M,\mathcal F)\to (M'',\mathcal F'')$ is a transversally harmonic map.
\end{coro}
\section{The first normal variational formula}
Let $(M,g,\mathcal F)$ be a foliated Riemannian manifold.
Let ${ vol}_L:M \to [0,\infty]$ be the volume function, of which ${
vol}_L (x)$ is the volume of the leaf  passing through $x\in
M$. It is trivial that $vol_L$ is a basic function. Then we have
the following.
\begin{lemma} \label{lem2-4} On a foliated Riemannian manifold $(M,\mathcal F)$, it holds that
\begin{align}
d_B vol_L + (vol_L) \kappa_B =0.
\end{align}
 \end{lemma}
\noindent{\bf Proof.} Let $\{ v_1, \cdots , v_p \}$ be linearly
independent vector fields of $\Gamma L$ such that
$vol_L = i(v_p) \cdots
i(v_1) \chi_{\mathcal F}$, where $\chi_{\mathcal F}$ is the
characteristic form of $\mathcal F$.  By the Rummler's formula [\ref{Tond}], $\varphi_0:=d\chi_{\mathcal F} + \kappa \wedge
\chi_{\mathcal F}$ satisfies  $i(v_p) \cdots i(v_1) \varphi_0 =0$.
Therefore we have
\begin{align*}
i(v_p) \cdots i(v_1) d \chi_{\mathcal F} =& -i(v_p) \cdots
i(v_1)( \kappa \wedge \chi_{\mathcal F}) \\
=& (-1)^{p+1} (vol_L) \kappa.
\end{align*}
On the other hand, a direct calculation gives
\begin{align*}
d (i(v_p) \cdots i(v_1) \chi_{\mathcal F})
&= (-1)^p i(v_p)  \cdots i(v_1) d \chi_{\mathcal
F}+\alpha(v_1,\cdots,v_p),
\end{align*}
where $\alpha(v_1,\cdots,v_p)=\sum_{j=1}^p (-1)^{p-j}i(v_p)\cdots
i(v_{j+1})\theta(v_j)\{i(v_{j-1}\cdots i(v_1)\chi_{\mathcal F}\}$.
Since $L$ is integrable, $\alpha(v_1,\cdots,v_p)\in L^*$ and so $\alpha(v_1,\cdots,v_p)=0$. Since $vol_L$ is a basic function, we have 
\begin{align*}
d_B vol_L &=d_B (i(v_p) \cdots i(v_1) \chi_{\mathcal F})\\
 &= (-1)^p i(v_p)\cdots i(v_1) d_B \chi_{\mathcal F}\\
&=-(vol_L) \kappa_B.
\end{align*}
So the proof is completed. $\Box$

Let $\Omega$ be a compact subset of $M$. Then the {\it transversal energy} of $\phi$ on $\Omega\subset
M$ is defined by
\begin{align}\label{eq2-4}
E_B(\phi;\Omega)={1\over 2}\int_{\Omega} | d_T \phi|^2 {1\over vol_L}
\mu_{M},
\end{align}
where $|d_T\phi|^2=\sum_a g_{Q'}(d_T\phi(E_a),d_T\phi(E_a))$ and $\mu_{M}$ is the volume element of $M$.

Let
$V\in\phi^{-1}Q'$. Obviously, $V$ may be considered as a vector
field on $Q'$ along $\phi$. Then there is a 1-parameter family of
foliated maps $\phi_t$ with $\phi_0=\phi$ and ${d\phi_t\over
dt}|_{t=0}=V$. The family $\{\phi_t\}$ is said to be a {\it foliated variation} of $\phi$ with the {\it normal variation vector field} $V$. Then we have the first normal variational formula(cf. [\ref{KW2}]).

\begin{thm} $(${\rm The first normal variational formula}$)$ Let $\phi:(M,\mathcal F)\to (M',\mathcal F')$
be a smooth foliated map, and all leaves of $\mathcal F$ be compact. Let $\{\phi_t\}$ be a smooth foliated variation of $\phi$ supported in a compact domain $\Omega$. Then
\begin{align}\label{eq2-5}
{d\over dt}E_B(\phi_t,\Omega)|_{t=0}=-\int_{\Omega} \langle
V,\tau_b(\phi)\rangle {1\over vol_L}\mu_{M},
\end{align}
where $V(x)={d\phi_t\over dt}|_{t=0}$ is the normal variation
vector field of $\{\phi_t\}$ and $\langle\cdot,\cdot\rangle$ is the pull-back metric on $\phi^{-1}Q'$.
\end{thm}
\noindent {\bf Proof.} Let $\Omega$ be a compact domain of $M$ and let $\{\phi_t\}$ be a foliated variation of $\phi$ supported in $\Omega$ with the normal variation vector field $V\in\phi^{-1}Q'$. Choose a local orthonormal basic frame
$\{E_a\}$ on $Q$ such that $(\nabla E_a)(x)=0$. Define
$\Phi:M \times (-\epsilon,\epsilon) \to M'$ by
$\Phi(x,t)=\phi_t (x)$ and set $E=\Phi^{-1}Q'$. Let $\nabla^\Phi$
denote the pull-back connection on $E$. Obviously,
$d_T\Phi(E_a)=d_T\phi_t (E_a)$ and $d\Phi({\partial\over\partial
t})={{d\phi_t}\over {dt}}$. Moreover, we have
$\nabla^\Phi_{\partial\over {\partial t}} {\partial\over{\partial
t}}=\nabla^\Phi_{\partial\over {\partial t}} E_a
=\nabla^\Phi_{E_a}{\partial\over{\partial t}}=0$. 
Hence we have
\begin{align*}
\frac{d}{dt} E_B(\phi_t,\Omega)
 &= \int_{\Omega}\sum_a
\langle\nabla^\Phi_{\partial\over{\partial t}} d_T \Phi(E_a), d_T
\Phi (E_a)\rangle
{1\over vol_L} \mu_{M} \\
&= \int_{\Omega} \sum_a\langle\nabla^\Phi_{E_a} d\Phi
(\frac{\partial}{\partial t}),
d_T \Phi (E_a)\rangle {1\over vol_L} \mu_{M}\\
&=\int_{\Omega} \sum_a\{E_a\langle \frac{d \phi_t}{dt}, d_T
\phi_t (E_a)\rangle - \langle \frac{d\phi_t}{dt}, \nabla_{E_a}^{\phi_t} d_T \phi_t (E_a)\rangle \} {1\over vol_L} \mu_{M}\\
&= \int_{\Omega}\sum_a E_a \{\langle \frac{d\phi_t}{dt},
d_T \phi_t (E_a)\rangle {1\over vol_L}\}\mu_{M}\\
&  - \int_{\Omega} \sum_a\langle \frac{d\phi_t}{dt},
d_T \phi_t (E_a)\rangle E_a({1\over vol_L}) \mu_{M} - \int_{\Omega} \langle \frac{d\phi_t}{d t}, \tau_b(\phi_t)\rangle
{1\over vol_L}\mu_{M}
\end{align*}
Now we define a normal vector field  $W_t$  by
\begin{align*}
W_t = {1\over vol_L}\sum_a\langle \frac{d\phi_t}{d t}, d_T \phi_t
(E_a)\rangle E_a.
\end{align*}
Then we have
\begin{align*}
 div_\nabla W_t  = \sum_a E_a
\{{1\over vol_L}\langle \frac{d\phi_t}{dt}, d_T \phi_t
(E_a)\rangle \}.
\end{align*}
 By the transversal divergence theorem (Theorem \ref{thm1-1}),  we have
 \begin{align*}
 \frac{d}{dt}E_B
(\phi_t,\Omega) &= \int_{\Omega} \{div_\nabla W_t- \langle
\frac{d\phi_t}{dt}, d_T \phi_t (d_B ({1\over vol_L}))\rangle\} \mu_{M}\\
&-\int_{\Omega} \langle \frac{d\phi_t}{dt}, \tau_b(\phi_t)\rangle
{1\over vol_L}\mu_{M}\\
&=\int_{\Omega}\langle {d\phi_t \over dt},d_T\phi_t((vol_L)\kappa_B +
d_B vol_L)\rangle{1\over vol_L^{2}}\mu_{M}\\
&-\int_{\Omega} \langle \frac{d\phi_t}{dt}, \tau_b(\phi_t)\rangle
{1\over vol_L}\mu_{M}.
\end{align*}
By Lemma \ref{lem2-4}, we have
\begin{align}
\frac{d}{dt}E_B (\phi_t,\Omega)=-\int_{\Omega} \langle \frac{d\phi_t}{dt},
\tau_b(\phi_t)\rangle {1\over vol_L}\mu_{M},
\end{align}
which proves (\ref{eq2-5}). $\Box$
\begin{coro} Let $\phi:(M,\mathcal F)\to (M',\mathcal F')$ be a smooth foliated map. Assume that all leaves of $\mathcal F$ are compact. Then $\phi$ is transversally harmonic if and only if $\phi$ is a critical point of the trasnversal energy of $\phi$ on any compact domain.
\end{coro}

\section{A generalized Weitzenb\"ock type formula and its applications}
Let $(M,g,\mathcal F)$ and $(M',g',\mathcal F')$ be two foliated Riemannian manifolds and let $\phi :(M,
\mathcal F) \rightarrow (M', \mathcal F')$ be a smooth foliated map. Note that $| d_T \phi |^2 \in
\Omega_B^0(\mathcal F)$ [\ref{KW1}]. Let $\Omega_B^r(E)=\Omega_B^r(\mathcal F)\otimes E$ be the space of $E$-valued basic $r$-forms, where $E=\phi^{-1}Q'$. We define $d_\nabla : \Omega_B^r(E)\to \Omega_B^{r+1}(E)$ by
\begin{align}
d_\nabla(\omega\otimes s)=(-1)^r\omega\wedge\nabla^\phi s + d_B\omega\otimes s
\end{align} 
for any $s\in E$ and $\omega\in\Omega_B^r(\mathcal F)$. 
Let $\delta_\nabla$ be a formal adjoint of $d_\nabla$. Then we have the following.
\begin{align}
d_\nabla = \sum_a \theta^a \wedge\tilde\nabla_{E_a}, \quad
\delta_\nabla = -\sum_a i(E_a) \tilde\nabla_{E_a} + i (\kappa_{B}^\sharp),
\end{align}
where $i(X)(\omega\otimes s)=i(X)\omega\otimes s$ for any $X\in TM$. Then the Laplacian $\Delta$ on $\Omega_B^*(E)$ is defined by 
\begin{align}
\Delta =d_\nabla \delta_\nabla +\delta_\nabla d_\nabla.
\end{align}
Moreover, the operators $A_X$ and $\theta(X)$ are extended to $\Omega_B^r(E)$ as follows:
\begin{align}
A_X(\omega\otimes s)&=A_X\omega\otimes s\\
\theta(X)(\omega\otimes s)&=\theta(X)\omega\otimes s+\omega\otimes\nabla_X^\phi s 
\end{align}
  for any $\omega\otimes s\in\Omega_B^r(E)$ and  $X\in TM$. Then $\theta(X)=d_\nabla i(X) +i(X)d_\nabla$ for any $X\in TM$. Hence $\Phi \in\Omega_B^*(E)$ if and only if $i(X)\Phi=0$ and $\theta(X)\Phi=0$ for all $ X\in \Gamma L$.
Then the  generalized Weitzenb\"ock type formula (\ref{eq1-12}) is extended to $\Omega_B^*(E)$ as follows:
\begin{align}\label{eq4-6}
\Delta \Phi = \tilde\nabla_{\rm tr}^*\tilde\nabla_{\rm tr} \Phi
+ F(\Phi) + A_{\kappa_{B}^\sharp} \Phi,\quad\forall \Phi\in\Omega_B^r(E),
\end{align}
where $F(\Phi)=\sum_{a,b=1}^{q}\theta^a\wedge i(E_b)\tilde R(E_b,E_a)\Phi$. Note that $d_T\phi\in\Omega_B^1(E)$. Then we have the following.
\begin{thm} Let $\phi:(M, g,\mathcal F) \to (M', g', \mathcal F')$ be a smooth foliated map. Then the generalized Weitzenb\"ock type formula is given by
\begin{align*}
\frac12\Delta_B| d_T \phi |^2 &= \langle\Delta d_T \phi, d_T \phi\rangle -
 |\tilde\nabla_{\rm tr} d_T \phi|^2 -
\langle A_{\kappa_{B}^\sharp}d_T \phi, d_T \phi\rangle  -\langle F(d_T\phi),d_T\phi\rangle,
\end{align*}
where 
\begin{align}\label{eq3-7-1}
\langle F(d_T\phi),d_T\phi\rangle&=\sum_a g_{Q'}(d_T \phi({\rm Ric^{Q}}(E_a)),d_T \phi(E_a))\\
&-\sum_{a,b}g_{Q'}( R^{Q'}(d_T \phi(E_b), d_T \phi(E_a))d_T \phi(E_a), d_T \phi(E_b)).\notag
\end{align}
\end{thm}
\noindent{\bf Proof}. Let $\{E_a\}(a=1,\cdots,q)$ be a local orthonormal basic frame such that 
at $x\in M$, $(\nabla E_a)_x=0$. Then, at $x$, we have from (\ref{eq1-7})
\begin{align}\label{eq5-8}
\frac12\Delta_B|d_T\phi|^2  &= \langle\tilde\nabla_{\rm tr}^* \tilde\nabla_{\rm tr} d_T \phi, d_T \phi\rangle - |\tilde\nabla_{\rm tr}d_T\phi|^2.
\end{align}
From (\ref{eq4-6}) and (\ref{eq5-8}), we have
\begin{align*}
\frac12\Delta_B |d_T\phi |^2 &=\langle\Delta d_T \phi, d_T \phi\rangle -
 |\tilde\nabla_{\rm tr} d_T \phi|^2 -
\langle A_{\kappa_{B}^\sharp}d_T \phi, d_T \phi\rangle  -\langle F(d_T\phi),d_T\phi\rangle.
\end{align*}
Now, we study $\langle F(d_T\phi),d_T\phi\rangle$. Let $\{V_\alpha\}(\alpha=1,\cdots,q')$ be a local orthonormal basic frame of $Q'$ and $\omega^\alpha$  be its dual coframe field. Let $f^\alpha=\phi^*\omega^\alpha$. Then $d_T\phi$ is expressed by
\begin{align}\label{eq4-7}
d_T\phi =\sum_{\alpha=1}^{q'}f^\alpha\otimes V_\alpha,
\end{align}
where $V_\alpha (x)\equiv V_\alpha(\phi(x))$. By direct calculation, we have
\begin{align}\label{eq4-8}
\tilde R(E_a,E_b)d_T\phi=\sum_{\alpha} R^{Q}(E_a,E_b)f^\alpha \otimes V_\alpha +\sum_\alpha f^\alpha \otimes R^E(E_a,E_b)V_\alpha,
\end{align}
where $R^E(E_a,E_b)V_\alpha =R^{Q'}(d_T\phi(E_a),d_T\phi(E_b))V_\alpha$.
From (\ref{eq4-8}), we have
\begin{align*}
\langle F(d_T\phi),d_T\phi\rangle&=\langle\sum_{a,b} \theta^a\wedge i(E_b)\tilde R(E_b,E_a)d_T\phi,d_T\phi\rangle\\
&=\sum_{a,b,\alpha,\beta}\langle \theta^a\wedge i(E_b)R^{Q}(E_b,E_a)f^\alpha\otimes V_\alpha,f^\beta\otimes V_\beta\rangle\\
&+\sum_{a,b,\alpha,\beta} g_{Q}(\theta^a\wedge i(E_b)f^\alpha,f^\beta)g_{Q}(R^E(E_a,E_b)V_\alpha,V_\beta).
\end{align*}
Note that $d_T\phi(E_a)=\sum_\alpha f^\alpha(E_a)V_\alpha$. Then we have
\begin{align}\label{eq4-10}
\sum_{a,b,\alpha}g_{Q}(\theta^a\wedge i(E_b)R^Q(E_b,E_a)f^\alpha,f^\alpha)=\sum_a g_{Q'}(d_T({\rm Ric}^{Q}(E_a)),d_T\phi(E_a)).
\end{align}
From (\ref{eq4-10}), we have
\begin{align*}
\langle F(d_T\phi),d_T\phi\rangle&=\sum_a g_{Q'}(d_T\phi({\rm Ric}^{Q}(E_a)),d_T\phi(E_a))\\
&+\sum_{a,b} g_{Q'}(R^{Q'}(d_T\phi(E_a),d_T\phi(E_b))d_T\phi(E_a),d_T\phi(E_b)),\notag
\end{align*}
which completes the proof. $\Box$

\noindent{\bf Remark}. (1) Let $\phi:(M,\mathcal F)\to (M',\mathcal F')$ be a smooth foliated  map. Then
\begin{align}
d_\nabla (d_T\phi)=0,\quad\delta_\nabla d_T\phi=-\tau_b (\phi) +i(\kappa_{B}^\sharp)d_T\phi.
\end{align}
(2) If a foliated map $\phi:(M, \mathcal F) \rightarrow (M',
\mathcal F')$ is  transversally harmonic, then
\begin{align}\label{eq4-12}
\Delta d_T \phi = d_\nabla i(\kappa_{B}^\sharp) d_T \phi.
\end{align}
\begin{coro} Let $\phi:(M,g,\mathcal F)\to (M',g',\mathcal F')$ be a transversally harmonic map. 
Then 
\begin{align}\label{eq4-13}
\frac12 \Delta_B |d_T\phi|^2 =-|\tilde\nabla_{\rm tr}d_T\phi|^2 -\langle F(d_T\phi),d_T\phi\rangle +\frac12 \kappa_{B}^\sharp (|d_T\phi|^2).
\end{align}
\end{coro}
{\bf Proof.} Since $d_\nabla(d_T\phi)=0$, we have
\begin{align}\label{eq4-14}
A_X d_T\phi =-\tilde\nabla_X d_T\phi + d_\nabla i(X) d_T\phi, \quad \forall X\in \Gamma Q.
\end{align}
Hence (\ref{eq4-13}) follows from (\ref{eq4-12}) and (\ref{eq4-14}). $\Box$

As applications of the generalized Weitzenb\"ock formula, we have the following theorems.
\begin{thm}
Let $(M,g,\mathcal F)$ be a compact foliated Riemannian manifold of nonnegative transversal Ricci curvature,  and  let $(M',g',\mathcal F')$ be a foliated Riemannian manifold of nonpositive transversal sectional curvature. If
$\phi:(M, \mathcal F) \rightarrow (M', \mathcal F')$ is
transversally harmonic,  then $\phi$ is transversally totally geodesic, i.e., $\nabla_{\rm tr}d_T\phi=0$.
Furthermore, \\
$(1)$ If the transversal Ricci curvature ${\rm Ric}^{Q}$ of $\mathcal F$ is
positive somewhere, then $\phi$ is transversally constant, i.e., $d\phi(TM)\subset L'$, equivalently, $\phi(M)\subset $ a leaf of $\mathcal F'$.
\\
$(2)$ If the transversal sectional curvature $K^{Q'}$ of $\mathcal F'$ is
negative at some point, then $\phi$ is either transversally constant or $\phi(M)$ is a transversally geodesic closed curve.
\end{thm}
\noindent{\bf Proof.} Let $\phi:(M,\mathcal F)\to (M',\mathcal F')$ be a transversally harmonic map. Then, from (\ref{eq4-13}), we have
\begin{align}\label{eq4-15}
\frac12 ( \Delta_B-\kappa_B^\sharp) |d_T\phi|^2 =-|\tilde\nabla_{\rm tr}d_T\phi|^2 -\langle F(d_T\phi),d_T\phi\rangle.
\end{align}
Since ${\rm Ric}^{Q}\geq 0 $ and $K^{Q'}\leq 0$, from (\ref{eq3-7-1}) we have
\begin{align}\label{eq4-16}
\langle F(d_T\phi),d_T\phi\rangle \geq 0.
\end{align}
Hence $(\Delta_B-\kappa_B^\sharp)|d_T\phi|^2\leq 0$. From Lemma \ref{lem2-3},  $|d_T\phi|$ is constant. Hence again, we have from (\ref{eq4-15})
\begin{align}
|\tilde\nabla_{\rm tr}d_T\phi|^2 +\langle F(d_T\phi),d_T\phi\rangle = 0.
\end{align}
Hence  $\tilde\nabla_{\rm tr}d_T\phi=0$ and
\begin{align}
&\sum_a g_{Q'}(d_T\phi({\rm Ric}^{Q}(E_a),d_T\phi(E_a))=0,\label{eq4-18}\\
&\sum_{a,b} g_{Q'}(R^{Q'}(d_T\phi(E_a),d_T\phi(E_b))d_T\phi(E_a),d_T\phi(E_b))=0.\label{eq4-19}
\end{align}
Therefore  $\phi$ is transversally totally geodesic. Moreover, from (\ref{eq4-18}), 
 if ${\rm Ric}^{Q}$ is positive at some point, then $d_T\phi=0$, i.e., $\phi$ is transversally constant, which proves (1). For the proof of (2), we assume that $K^{Q'}<0$. From (\ref{eq4-19}),  the rank of $d_T\phi <2$. Hence the rank of $d_T\phi$ is zero or one everywhere. If the rank of $d_T\phi$ is zero, then $\phi$ is transversally constant. If the rank of $d_T\phi$ is one, then $\phi(M)$ is closed transversally geodesic. $\Box$
 
 Next, we extend Theorem 5.3 under the weaker transversal sectional curvature of $\mathcal F'$. Let $rank_T(\phi)$ be the rank of $d_T\phi$.
 \begin{thm} Let $(M,g,\mathcal F)$ be a compact foliated Riemannian manifold and let $(M',g',\mathcal F')$ be a foliated Riemannian manifold. Assume that $Ric^{Q} \geq\lambda$ {\rm id.}  and $K^{Q'}\leq\mu$ for any positive constants $\lambda$ and $\mu$. Let $\phi : (M,\mathcal F) \rightarrow (M', \mathcal F')$ be a transversally
harmonic map with {\rm max}$\{rank_T (\phi)\}\leq C$, where $C\geq 2$ is constant. If $|d_T\phi|^2 \leq \frac{\lambda C}{\mu (C-1)}$, then $\phi$ is transversally constant or $\phi$ is  transversally geodesic. In
particular, if $|d_T\phi|^2 \leq \frac{\lambda}{\mu}$, then $\phi$ is  transversally constant.
\end{thm}
{\bf Proof.} Let $\{E_a\}$ be a local orthonormal basic frame of $Q$. From (\ref{eq3-7-1}), we have
\begin{align}\label{eq4-20}
\langle F(d_T\phi),d_T\phi\rangle &= \sum_a g_{Q'}(d_T\phi({\rm Ric}^{Q}(E_a)),d_T\phi(E_a))\\
 &-\sum_{a,b}\{|d_T\phi(E_a)|^2|d_T\phi(E_b)|^2-g_{Q'}(d_T\phi(E_a),d_T\phi(E_b))^2\}K_{ab}^{Q'},\notag
 \end{align}
 where $K_{ab}^{Q'}= K^{Q'}(d_T\phi(E_a),d_T\phi(E_b))$ is the transversal sectional curvature spanned by $d_T\phi(E_a)$ and $d_T\phi(E_b)$.
Let $rank_T(\phi)=r\leq C$. Now, we choose a local orthonormal basic frame $\{E_a\}$ such that 
$g_{Q'}(d_T(E_a),d_T(E_b))|_x =\lambda_a\delta_{ab}$ and $\lambda_1\geq\lambda_2\geq \cdots\geq\lambda_r>0$.
Then, from (\ref{eq4-13}) and (\ref{eq4-20}), we have
\begin{align*}
\frac12\Delta_B |d_T\phi|^2 &= -|\tilde\nabla_{\rm tr}d_T\phi|^2 +\frac12 \kappa_{B}^\sharp (|d_T\phi|^2)
-\sum_a g_{Q'}(d_T\phi({\rm Ric}^{Q}(E_a)),d_T\phi(E_a))\\
 &+\sum_{a,b}\{|d_T\phi(E_a)|^2|d_T\phi(E_b)|^2-g_{Q'}(d_T\phi(E_a),d_T\phi(E_b))^2\}K_{ab}^{Q'}\\
&\leq -|\tilde\nabla_{\rm tr} d_T \phi|^2  +\frac12 \kappa_{B}^\sharp (|d_T\phi|^2)- \lambda |d_T\phi|^2 + \mu \{|d_T\phi|^4
- \sum_{a=1}^r \lambda_a^2 \}
\end{align*}
Using the Schwarz's inequality, we have
\begin{align}\label{eq4-21}
|d_T\phi|^4 = \sum_{a,b}^r \lambda_a \lambda_b
\leq \frac12
\sum_{a,b=1}^r (\lambda_a^2 + \lambda_b^2)
= r \sum_{a=1}^r \lambda_a^2 \leq C \sum_{a=1}^r
\lambda_a^2.
\end{align}
From (\ref{eq4-21}), we have
\begin{align}\label{eq4-22}
 |d_T\phi|^4 - \sum_{a=1}^r \lambda_a^2 \leq \frac{C-1}C |d_T\phi|^4.
\end{align}
From (\ref{eq4-22}), we have
\begin{align}\label{eq4-23}
\frac12\Delta_B |d_T\phi|^2 &\leq  -|\tilde\nabla_{\rm tr} d_T \phi|^2  +\frac12 \kappa_{B}^\sharp (|d_T\phi|^2)- |d_T\phi|^2\{\lambda - {{(C-1)\mu}\over C}|d_T\phi|^2\}\notag\\
&\leq \frac12 \kappa_{B}^\sharp (|d_T\phi|^2).
\end{align}
Hence, from Lemma \ref{lem2-3},  $|d_T\phi|$ is constant and then
\begin{align}
|\tilde\nabla_{\rm tr} d_T \phi|^2+ |d_T\phi|^2\{\lambda - {{(C-1)\mu}\over C}|d_T\phi|^2\}=0.
\end{align}
Therefore $\tilde\nabla_{\rm tr}d_T\phi=0$ and $ |d_T\phi|^2\{\lambda - {{(C-1)\mu}\over C}|d_T\phi|^2\}=0$. Hence $\phi$ is transversally totally geodesic.
If $d_T\phi=0$, then $\phi$ is transversally constant. If $d_T\phi\ne 0$, then $|d_T\phi|^2={\lambda C \over \mu(C-1)}$ and $\phi$ is transversally totally geodesic. In particular,
if $|d_T\phi|^2 \leq {\lambda\over\mu}$, then $\phi$ is transversally constant. $\Box$
 
\noindent{\bf Remark.} For the point foliation, Theorem 5.3 and Theorem 5.4 are found in [\ref{ES}] and [\ref{Se}], respectively.  

\vskip 0.3cm
\noindent{\bf Acknowledgements} This research was supported by the Basic Science Research Program
        through the National Research Foundation of Korea(NRF) funded
        by the Ministry of Education, Science and Technology(2010-0021005).

\noindent Min Joo Jung

\noindent Department of Mathematics, Jeju National University,
Jeju 690-756, Korea

\noindent{\it E-mail address} :niver486@jejunu.ac.kr

\vskip0.3cm

\noindent Seoung Dal Jung

\noindent Department of Mathematics and Research Institute for Basic Sciences, Jeju National University,
Jeju 690-756, Korea

\noindent {\it E-mail address} : sdjung@jejunu.ac.kr

\end{document}